\documentclass[12]{amsart}

\usepackage{color}

\newtheorem{define}{Definition}
\newtheorem{proposition}[define]{Proposition}
\newtheorem{conjecture}[define]{Conjecture}

\newtheorem{theorem}[define]{Theorem}

\newcommand{\bd}{\begin{define} \rm}
	\newcommand{\ed}{\end{define}}

\newfont{\Bb}{msbm10 scaled\magstep1}

\begin{document}
\title[Hedetniemi's and Stahl's conjectures and the Poljak-R\"odl function]{A note 
on Hedetniemi's conjecture, Stahl's conjecture and the Poljak-R\"odl function}
\author{Claude Tardif}
\address{Royal Military College of Canada} 
\email{Claude.Tardif@rmc.ca} 
\author{Xuding Zhu}
\address{Zhejiang Normal University} 	
		\email{xdzhu@zjnu.edu.cn} 
		\thanks{Xuding Zhu's research is supported by grant mumbers NSFC 11971438, ZJNSF LD19A010001 and 111 project of the Ministry of Education of China.}

	\maketitle
	\begin{abstract}
		We prove that $\min\{\chi(G), \chi(H)\} - \chi(G\times H)$ can be arbitrarily large,
		and that if Stahl's conjecture on the multichromatic number of Kneser graphs holds, 
		then we can have $\chi(G\times H)/\min\{\chi(G), \chi(H)\} \leq 1/2 + \epsilon$
		for large values of $\min\{\chi(G), \chi(H)\}$.
	\end{abstract}
	
	\section{Introduction}
	
	The {\em categorical product} $G \times H$ of graphs $G$ and $H$ has vertex set $V(G \times 
	H)=\{(x,y): x\in V(G), y \in V(H)\}$, in which two vertices $(x,y)$ and $(x',y')$ are adjacent if and only if 
	$xx' \in E(G)$ and $yy' \in E(H)$.
	A proper colouring $\phi$ of $G$ can be lifted to a proper colouring 
	$\Phi$ of $G \times H$ defined as $\Phi(x,y)=\phi(x)$. So $\chi(G \times H) \le  \chi(G) $, and similarly 
	$\chi(G \times H) \le \chi(H)$.  
	Hedetniemi  conjectured in 1966  that $\chi(G \times H) = \min \{\chi(G),\chi(H)\}$ for all finite graphs 
	$G$ and $H$ \cite{Hedetniemi}. 
	The conjecture received a lot of   attention \cite{Klav, Sauer,Tardif,Zhu} and remained open for more 
	than half century. It is known that $\chi(G\times H) = \min\{\chi(G), \chi(H)\}$ 
	whenever $\min\{\chi(G), \chi(H)\} \le 4$ \cite{ES} and that the fractional version is true, i.e., for any 
	graphs $G$ and $H$, $\chi_f(G \times H) = \min \{\chi_f(G),\chi_f(H)\}$ \cite{Zhu-frac}.
	However,  Shitov  refuted this conjecture recently \cite{Shitov}. Yet, some problems concerning the 
	chromatic number of product graphs remain open. 
	
	The {\em Poljak-R\"odl function} $f: \mbox{\Bb N} \rightarrow \mbox{\Bb N}$ is defined by
	$$f(n) = \min\{\chi(G \times H): \chi(G), \chi(H) \ge n\}.$$
	%As observed above,  $f(n) \le n$ for all $n$.
	Hedetniemi's conjecture is equivalent to the statement that $f(n)=n$ for all $n$. 
	Shitov proved that for sufficiently large $n$, $f(n) \le n-1$. 
	Still, very little is known about the behavior of the function $f(n)$. In particular, it is unknown whether 
	$f(n)$ is bounded by a constant. However it is known that if $f(n)$ 
	is bounded by a constant, then $f(n) \le 9$ for all $n$ (see \cite{Sauer, Zhu}).
	In this note, we prove the following facts.
	
	\begin{proposition} \label{aboprf}
		
		\begin{itemize}
			\item[] 
			\item[(i)] $\lim_{n \rightarrow \infty} (n-f(n)) = \infty$,
			\item[(ii)] If Stahl's conjecture on the multichromatic number of Kneser graphs 
			\cite{Stahl} holds, then $\limsup_{n \rightarrow \infty} f(n)/n \leq 1/2$.
		\end{itemize}
	\end{proposition}
	
	Proposition~\ref{aboprf} (i) will be proved in Section~\ref{de}. Proposition~\ref{aboprf} (ii)
	will be proved in Section~\ref{sc}, where a presentation of Stahl's conjecture 
	is also given.

	\section{Discussion and extensions of Shitov's results} \label{de}
	
	For a positive integer $c$, the {\em exponential graph} $K_c^H$ has vertices 
	all the mappings $f: V(H) \to \{1,2,\ldots, c\}$, in which $f, g$ are adjacent
	in $K_c^H$ if $f(u) \ne g(v)$ for every edge $e=uv$ of $H$. It is well known 
	and easy to verify that 
	$\Phi(v,f)=f(v)$ is a proper $c$-colouring of $H \times K_c^H$.
	Thus the way to find counterexamples to Hedetniemi's conjecture is to find
	an integer $c$ and a graph $H$ such that both $H$ and $K_c^H$ have chromatic 
	number larger than $c$. 

The {\em lexicographic product} $G[H]$ of $G$ and $H$ is the graph with vertex set 
$V(G[H])=\{(x,y): x\in V(G), y \in V(H)\}$, in which two vertices $(x,y)$ and 
$(x',y')$ are adjacent if and only if $xx' \in E(G)$, or $x = x'$ and 
$yy' \in E(H)$. 

Shitov's construction of counterexamples to Hedetniemi's conjecture is based on the following result.
	\begin{theorem}[\cite{Shitov}, Claim 3] \label{claim3}
		For any graph $G$ with girth at least six, 
		for all but finitely many values of $q$, we have 
		$\chi \left  ( K_c^{G[K_q]} \right ) \geq c+1$, with $c = \lceil 3.1 q \rceil$\footnote{Technically, 
Shitov refers to the ``strong product'' rather than the lexicographic product of graphs, but with $K_q$ as 
a second factor, the strong product coincides with the lexicographic product (see \cite{IK}).}.
	\end{theorem}

Finding such a lower bound on chromatic numbers of some exponential graphs was the key part of
Shitov's refutation of Hedetniemi's conjecture. Finding lexicographic products $G[K_q]$ 
with $\chi(G[K_q]) > c$ is standard theory. Indeed the fractional chromatic number 
$\chi_f(H)$ of a graph $H$ is a standard lower bound for its chromatic number, and it is well known 
that $\chi_f(G[H]) = \chi_f(G) \chi_f(H)$ (see \cite{GZ}).  Erd\H{o}s' classic probabilistic proof \cite{Erdos}
shows that there are graphs with girth at least $6$ and fractional chromatic number at least $3.1$.
For such a graph $G$, we have $\chi(G[K_q]) \geq \lceil \chi_f(G[K_q]) \rceil = \lceil\chi_f(G) \cdot q \rceil
\geq \lceil 3.1 q \rceil$,
and by Theorem \ref{claim3}, this yields a counterexample to Hedetniemi's conjecture.

Remarkably, replacing the condition $\chi_f(G) \geq 3.1$ by $\chi_f(G) \geq B$ for $B \gg 3.1$ readily gives 
counterexamples to Hedetniemi's conjecture where the chromatic number of at least one factor is arbitrarily 
larger than the chromatic number of the product. Also, the proof of Theorem~\ref{claim3} only uses 
a small subgraph of $K_c^{G[K_q]}$. Therefore it is possible that Shitov's construction already gives 
examples that show that $\lim_{n \rightarrow \infty} f(n)/n = 0$. On the other hand, since $\chi_f(G[K_q]) > c$,
the fractional version of Hedetniemi's conjecture \cite{Zhu-frac} implies that $\chi_f(K_c^{G[K_q]}) = c$.
Thus it is also reasonnable to think that  $\chi(K_c^{G[K_q]})/c$ may be bounded,
and that the identity $\lim_{n \rightarrow \infty} f(n)/n = 0$, if true, can only be witnessed by a different construction.

	\smallskip \noindent \begin{proof}[{Proof of Proposition~\ref{aboprf} (i).}]
	Fix a positive integer $d$. We shall prove that if $n$ is sufficiently large, then 
	$f(n+d) \le n$.
	Let $G_d$ be a graph with girth at least 6 and fractional 
	chromatic number at least $8d$. Then by Theorem~\ref{claim3}, for sufficiently large $q$ and 
	$c = \lceil 3.1 q \rceil$, we have $\chi \left  ( K_c^{G_d[K_q]} \right ) \geq c+1$
	while $\chi(G_d[K_q]) \geq 2cd$. Now consider the graph
	$K_{cd}^{G_d[K_q]}$. For $i=0,1,\ldots, d-1$, let $Q_i$ be the subgraph of 
		$K_{cd}^{G_d[K_q]}$   induced by the functions
	with image in $\{ic+1,ic+2, \ldots, ic+ c\}$. Each $Q_i$ is isomorphic to $K_c^{G_d[K_q]}$
	and hence at least $c+1$ colours are needed for each copy. For $i \ne j$, 
	each function in $Q_i$ is adjacent to each function in $Q_j$.
	Hence, $\chi \left  ( K_{dc}^{G_d[K_q]} \right ) \geq d(c+ 1)$. As $\chi(K_{dc}^{G_d[K_q]}) = dc$ and 
		$\chi(G_d[K_q]) \ge 2cd \ge cd+d$, it follows
	that $f(dc + d) \leq dc$. 

Thus for every $d$ there exist infinitely many values of $n$ (of the form $dc + d$) such that 
$n - f(n) \geq d$. It only remains to show that the gap between $n$ and $f(n)$ will not close
while going from one value of $c$ to the next. Note that $c = \lceil 3.1 q \rceil$, where $q$ 
is any value above a fixed threshold, and $\lceil 3.1 (q+1) \rceil - \lceil 3.1 q \rceil \leq 4$.
Thus it suffices to examine the values $n = dc + d + i$ where $i \leq 4d$, and we can suppose that $c \geq 5$.
The graph $K_{cd + i}^{G_d[K_q]}$ contains a copy of 
		$K_{cd}^{G_d[K_q]}$   induced by the functions
		with image in $\{1,2, \ldots, cd\}$. For $j=cd+1, cd+2, \ldots,
		cd+i$, the constant functions $g_j$ with image $j$ are pairwise adjacent and each is adjacent to 
		all the functions in  $K_{cd}^{G_d[K_q]}$. Hence $\chi(K_{cd + i}^{G_d[K_q]}) \ge \chi(K_{cd}^{G_d[K_q]})+i \ge cd+d+i$. 
		For $i \leq (c-1)d$, we also have $\chi(G_d[K_q]) \geq cd+d+i$, so that $f(cd+d+i) \leq cd+i$.
		Altogether, the inequality $f(n+d) \leq n$ is established for all 
		but finitely many values of $n$. Thus, $\lim_{n\rightarrow \infty} n - f(n) = \infty$.
	\end{proof}

The gap between $n$ and $f(n)$ proved in this section depends on the minimum number $p$ of vertices 
	of a   girth $6$ graph  with fractional chromatic number at least $8d$. The best known  upper bound for 
	$p$ to our knowledge  is  $p = O((d \log d)^4)$, which follows from a result of Krivelevich \cite{Kri}. Using this result, one can  show that
	for any   $\epsilon > 0$, there is a constant $a$ such that for sufficiently large $n$, 
	$f(n) \le n- a (\log n)^{1/4 - \epsilon}$.
Very recently, He and Wigderson \cite{HW} proved that for some $\epsilon \simeq 10^{-9}$,   $f(n) < (1 - \epsilon)n$
for sufficiently large $n$.
The examples are again cases of Shitov's construction.

	\section{Stahl's conjecture} \label{sc}
	
	In the proof of Proposition 1(i), based on the fact that
	$\chi(K_c^{G_d[K_q]})\geq c+1$, we have shown that 
	$\chi(K_{cd}^{G_d[K_q]})\geq cd+d$. In this section, we show that if a
	special case of a conjecture of Stahl on the 
	multichromatic number of Kneser graphs is true, then $\chi(K_{cd}^{G_d[K_q]})$
	is much larger.
	
	Consider a proper colouring $\phi$ of the graph $K_{cd}^{G_d[K_q]}$
	with $x$ colours. Let $A$ be a subset of
	$\{1, \ldots, cd\}$ of cardinality $c$. Let $R_A$ be the subgraph of
	$K_{cd}^{G_d[K_q]}$ induced by the functions with image contained in $A$.
	Then $R_A$ is isomorphic to $K_{c}^{G_d[K_q]}$, so $\phi$ uses at least $c+1$
	colours on $R_A$. Let $\psi(A) \subseteq \{1,\ldots,x\}$ be a subset of
	exactly $c+1$ colours used by $\phi$ on $R_A$. We have $\psi(A)$ disjoint 
	from $\psi(B)$ whenever $A$ is disjoint from $B$, because $R_A$ is totally
	joined to $R_B$ in $K_{cd}^{G_d[K_q]}$. This property can be formulated 
	in terms of homomorphisms of Kneser graphs. Recall
	that the vertices of the Kneser graph $K(m,n)$ are the 
	$n$-subsets of $\{1, \ldots, m\}$, and two of these are joined by 
	an edge whenever they are disjoint. Thus the colouring 
	$\phi: K_{cd}^{G_d[K_q]} \rightarrow K_x$ induces a 
	homomorphism $\psi: K(cd,c) \rightarrow K(x,c+1)$.
	The question is how large does $x$ need to be for such a homomorphism to exist.
	
	Stahl's conjecture deals with the latter question. For an integer $n$,
	the {\em $n$-th multichromatic number} $\chi_n(H)$ of a graph $H$
	is the least integer $m$ such that $H$ admits a homomorphism to $K(m,n)$.
	In particular $\chi_1(H) = \chi(H)$. 
	Lov\'asz \cite{Lovasz} proved that 
	$\chi_1(K(m,n)) = \chi(K(m,n)) = m - 2n + 2$.
	Stahl~\cite{Stahl} investigated the general multichromatic numbers of Kneser graphs,
	and observed the following.
	\begin{itemize}
		\item[{\bf a.}] For $1 \leq k \leq n$, $\chi_k(K(m,n)) = m - 2(n-k)$,
		\item[{\bf b.}] $\chi_{kn}(K(m,n)) = km$,
		\item[{\bf c.}] $\chi_{k + k'}(K(m,n)) \leq \chi_k(K(m,n)) + \chi_{k'}(K(m,n))$.
	\end{itemize}
	Based on this he conjectured the following.
	\begin{conjecture}[\cite{Stahl}] If $k = an + b$, $a \geq 1, 0 \leq b \leq n-1$,
		then for $m \geq 2n$,
		$$\chi_k(K(m,n)) = \chi_{an}(K(m,n)) + \chi_b(K(m,n)) = (a+1)m - 2(n-b).$$
	\end{conjecture}
	
	\smallskip \noindent \begin{proof}[{Proof of Proposition~\ref{aboprf} (ii).}] 
	For a fixed $d$, let $G_d$ have girth at least $6$ and fractional chromatic number at least $8d$.
	For any $q$ above a given threshold $q_d$ and for $c = \lceil 3.1 q \rceil$,
	we have $\chi(G_d[K_q]) \geq 2cd$ and $\chi \left ( K_{cd}^{G_d[K_q]} \right )
	\geq \chi_{c+1}(K(cd,c))$, as explained in the first three paragraphs of this section.
	If Stahl's conjecture holds, then 
	$\chi \left ( K_{cd}^{G_d[K_q]} \right ) \geq 2cd - 2c + 2$.
	Since $f$ is monotonic, this gives $f(2cd - 2c + 2) \leq cd$.
	Therefore 
	$$n \in [2(c-4)d - 2(c-4) + 2, 2cd - 2c + 2] \mbox{\rm \  implies }
	\frac{f(n)}{n} \leq  \frac{cd}{2(c-4)d -2(c-4) + 2}.$$
	The intervals $[2(c-4)d - 2(c-4) + 2, 2cd - 2c + 2], c \in \mbox{\Bb N}$ cover all
	but a finite part of $\mbox{\Bb N}$. Hence 
	$$\limsup \frac{f(n)}{n} \leq \lim_{c \rightarrow \infty} \frac{cd}{2(c-4)d -2(c-4) + 2} 
	= \frac{d}{2d-2}.$$ 
	Since this holds for arbitrarily large $d$,  $\limsup \frac{f(n)}{n} \leq \frac 12$. \end{proof}

	\smallskip \noindent 
	{\bf Acknowledgements.} We thank Yaroslav Shitov for many helpful comments.

\end{document}